\documentclass{amsart}

\usepackage[utf8]{inputenc}

\usepackage{url}

\usepackage{units}

\usepackage[all]{xy}
\usepackage{pb-diagram, pb-xy}

\usepackage{mathtools}

\usepackage{amssymb}
\usepackage{amsthm}

\usepackage{mathrsfs}

\newtheorem{thm}{Theorem}[section]
\newtheorem{lem}[thm]{Lemma}
\newtheorem{cor}[thm]{Corollary}
\newtheorem{prop}[thm]{Proposition}

\theoremstyle{definition}
\newtheorem{ex}[thm]{Example}
\newtheorem{dfn}[thm]{Definition}

\newcommand{\bbC}{{\mathbb C}}

\newcommand{\bbN}{{\mathbb N}}

\newcommand{\bbQ}{{\mathbb Q}}

\newcommand{\bbZ}{{\mathbb Z}}

\newcommand{\frakS}{{\mathfrak S}}

\newcommand{\Supp}{\mathrm{Supp}}

\newcommand{\Pic}{{\mathit{Pic}}}

\DeclareMathOperator{\rk}{\mathit{rk}}

\DeclareMathOperator{\Dbc}{D^\mathit{b}_\mathit{c}}
\DeclareMathOperator{\Perv}{P}

\DeclareMathOperator{\Hom}{Hom}

\DeclareMathOperator{\Aut}{Aut}

\DeclareMathOperator{\rmD}{D}
\DeclareMathOperator{\rmF}{F}
\DeclareMathOperator{\rmH}{H}
\DeclareMathOperator{\rmK}{K}
\DeclareMathOperator{\rmN}{N}

\DeclareMathOperator{\Rep}{Rep}

\DeclareMathOperator{\Pervc}{{P_{\hspace*{-0.2em} \mathit{sc}}}}
\newcommand{\conv}{\circ}

\newcommand{\chiloc}{\chi_{loc}}

\newcommand{\tdual}{\hat{\theta}}

\newcommand{\Adual}{\hat{A}}

\newcommand{\calF}{{\mathcal F}}
\newcommand{\calH}{{\mathcal H}}

\newcommand{\calFtilde}{\tilde{\calF}}

\newcommand{\scrC}{{\mathscr C}}
\newcommand{\scrS}{{\mathscr S}}

\newcommand{\ftilde}{{\tilde{f}}}
\newcommand{\rtilde}{{\tilde{r}}}

\newcommand{\Ztilde}{{\tilde{Z}}}

\numberwithin{equation}{thm}

\begin{document}

\title[Characteristic classes and Hilbert-Poincar\'e series]{Characteristic classes and Hilbert-Poincar\'e series for perverse sheaves on abelian varieties} 
\author{Thomas Kr\"amer}
\address{Mathematisches Institut\\ Ruprecht-Karls-Universit\"at Heidelberg\\ Im Neuenheimer Feld 205, D-69120 Heidelberg, Germany}
\email{tkraemer@mathi.uni-heidelberg.de}


\begin{abstract}
The convolution powers of a perverse sheaf on an abelian variety define an interesting family of 
branched local systems whose geometry is still poorly understood. 
We show that the generating series for their generic rank is a rational function of a very simple shape and that a similar result holds for the symmetric convolution powers. We also give formulae for other Schur functors in terms of characteristic classes on the dual abelian variety, and as an example we discuss the case of Prym-Tjurin varieties.
\end{abstract}

\maketitle

\thispagestyle{empty}

\section{Introduction}

Tannakian categories arise naturally in many areas of algebraic geometry. For instance, Gabber and Loeser have 
 studied convolutions of perverse sheaves on tori via certain Tannakian quotient categories~\cite{GaL} which have also been used recently by Katz in his work on the Mellin transform over finite fields~\cite{KatzSatoTate}. In what follows we will be concerned with complex abelian varieties, where similar results have been obtained in relation with a generic vanishing theorem for perverse sheaves~\cite{KrWVanishing}. The Tannakian formalism produces many new objects from a given one by considering convolution powers and their subquotients such as symmetric powers; the arising perverse sheaves encode nontrivial information on the geometry of moduli spaces and Albanese morphisms~\cite{KrCubic}~\cite{KrE6}~\cite{KrWSchottky}~\cite{WeBN}~\cite{WeTorelli}, but they are hard to approach explicitly. 
The goal of this note is to study their generic rank. It turns out that the generating series for the symmetric  powers is a rational function of a very simple shape reminiscent of the Hilbert-Poincar\'e series of a graded module.

\medskip

Let $A$ be a complex abelian variety, and denote by $\Dbc(A)=\Dbc(\bbC_A)$ the bounded derived category of the category of constructible sheaves of complex vector spaces on $A$. The group law~$a: A\times A \to A$ of the abelian variety defines a convolution product
\[
K*L \;=\; Ra_*(K\boxtimes L) \quad \textnormal{for} \quad K, L \in \Dbc(A),
\]
and with respect to this convolution product the derived category becomes a rigid symmetric monoidal category. The abelian subcategory $\Perv(A) \subset \Dbc(A)$ of perverse sheaves is not stable under convolution, but one gets a Tannakian quotient category as follows. By the vanishing theorem of~\cite{KrWVanishing}~\cite{SchnellHolonomic} the Euler characteristic of perverse sheaves is nonnegative:
\[ 
 \chi(P) \;\;=\;\; \sum_{n\in \bbZ} \; (-1)^n \dim_\bbC \rmH^n(A, P) \;\; \geq \;\; 0
 \quad \textnormal{for all} \quad P\in \Perv(A).
\] 
We say that a complex $K\in \Dbc(A)$ is {\em negligible} if the Euler characteristic of each of its perverse cohomology sheaves vanishes. The negligible complexes form a thick subcategory $\rmN(A)\subset \Dbc(A)$, and the quotient category $\Perv(A)/(\Perv(A)\cap \rmN(A))$ is a limit of Tannakian categories with a tensor product induced from the convolution product~\cite{KrWVanishing}. If $\Pervc(A)$ denotes the category of {\em semisimple} perverse sheaves which are {\em clean} in the sense that they have no negligible subobjects, then for $P, Q \in \Pervc(A)$ we have a unique decomposition
\[
 P*Q \;=\; P\conv Q \;\oplus\; P\bullet Q
 \quad \textnormal{with} \quad
 \begin{cases}
  P\conv Q \;\in\; \Pervc(A), \\
  P \bullet Q \;\in\; \rmN(A),
 \end{cases}
\]
and $\Pervc(A)$ equipped with the tensor product $\conv$ is an inductive limit of Tannakian categories which lift the above Tannakian quotient categories. A similar lift exists in the non-semisimple case~\cite[sect.~5]{KraemerSemiabelian}, but working with semisimple perverse sheaves has the advantage that for these the decomposition theorem and the relative hard Lefschetz theorem holds by Kashiwara's conjecture, see~\cite{DrinfeldKashiwara}~\cite{BKdeJong}~\cite{GaitsgoryDeJong} or~\cite{MochizukiAsymptotic}~\cite{SabbahPolarizable}.

\medskip

Now the objects we are interested in arise in the above Tannakian framework as follows. For $P\in \Pervc(A)$ the symmetric group $\frakS_n$ acts on $P^{n} = P \conv \cdots \conv  P\in \Pervc(A)$ by permutation of the factors. Up to isomorphism the irreducible finite-dimensional rational representations $V_\alpha  \in \Rep_\bbQ(\frakS_n)$ are parametrized by the partitions $\alpha$ of $n$, and by~\cite[sect.~1.4]{DelCT} 

\[
 P^{n} \;\; \cong  \bigoplus_{\deg \alpha = n}  S^{\alpha}(P) \otimes_\bbQ V_\alpha
 \quad \textnormal{for the Schur functors} \quad
 S^{\alpha}: \; \Pervc(A)\;  \rightarrow \; \Pervc(A).
\]
For the partition $\alpha = (n)$ we get the symmetric convolution powers~$S^n(P)$ which are the invariants under the symmetric group. 
Notice that if the perverse sheaf $P$ underlies a mixed Hodge module or a polarized twistor module, then by the axioms in~\cite[sect.~5]{KrWVanishing} the same will also hold for all the perverse sheaves $S^\alpha(P)$ since the corresponding categories are $\bbQ$-linear pseudoabelian.

\medskip

Unfortunately, very little is known about the perverse sheaves $P^{n}$ and $S^{\alpha}(P)$ in general. Put $g=\dim A >0$, and let $U_n\subseteq A$ be an open dense subset over which the cohomology sheaves
\[
 \calF_n(P) \;=\; \calH^{-g}(P^{n} )|_{U_n} \quad \; \textnormal{and} \; \quad \calFtilde_\alpha \;=\; \calH^{-g}(S^{\alpha}(P))|_{U_n}
\]
are locally constant. The underlying monodromy representations are unknown even in the simplest nontrivial cases. The goal of this note is to determine the generic ranks
$r_P(n) = \rk(\calF_n(P))$ and $\rtilde_P(\alpha) = \rk(\calFtilde_\alpha(P))$. We will see that the generating series for the convolution powers and for the symmetric convolution powers
\[
 Z_P(t) \;=\; \sum_{n=1}^\infty \; r_P(n) \cdot t^n
 \quad \; \textnormal{and} \; \quad
 \tilde{Z}_P(t) \;=\; \sum_{n=1}^\infty \; \tilde{r}_P(n) \cdot t^n 
\]
are rational functions of a very simple shape. To control the negligible summands we will assume that the  {\em spectrum}~$\scrS(P)$ as defined in~\cite{KrWVanishing} is finite. Recall that~$\scrS(P)$ is the set of all characters 
$$\varphi \;\in\; \Pi(A) \;=\; \Hom_\bbZ(\pi_1(A, 0), \bbC^*)$$ 
such that the corresponding rank one local system $L_\varphi$ satisfies $\rmH^i(A, P\otimes_\bbC L_\varphi) \neq 0$ for some $i\neq 0$. The generic vanishing theorem of loc.~cit.~says that this spectrum is always a finite union of translates of algebraic subtori $\Pi(B) \subset \Pi(A)$ defined by proper abelian quotient varieties $A\twoheadrightarrow B$, and our finiteness assumption means that only $B=\{0\}$ occurs. This assumption is automatically satisfied if~$A$ is a simple abelian variety but also in many other applications, like for the perverse intersection cohomology sheaves on curves in lemma~\ref{lem:spectrum} below or on ample divisors with only isolated singularities~\cite[lemma 8.1]{KrE6}. This being said, we obtain the following result for the generating series for the generic ranks from above.

\begin{thm} \label{thm:polynomial}
Let $P\in \Pervc(A)$ be a semisimple clean perverse sheaf with a finite spectrum and with Euler characteristic $\chi=\chi(P)$. Then 
\medskip
\begin{align} \nonumber
 Z_P(t) & \;=\; \frac{t f_P(t)}{(1-\chi t )^{g+1}} \quad
 \textnormal{\em for some  $f_P(t) \in \bbZ[t]$  of degree $\deg_t f_P(t) \leq g-1$}, \\[0.5em] \nonumber 
 \Ztilde_P(t) & \;=\; \frac{t\ftilde_P(t)}{(1-t)^{2g+\chi}} \quad
 \textnormal{\em for some  $\ftilde_P(t) \in \bbZ[t]$ of degree $\deg_t \ftilde_P(t) \leq 2g-2$}, 
\end{align} 
\smallskip \\
\noindent
and we have the functional equation \medskip
\[
 \ftilde_P(t) \;=\; t^{2g-2}\, \ftilde_P(1/t).
\]
\end{thm}

\medskip

Notice that since we have introduced an extra factor $t$ in the nominator of our rational functions, the constant term $f_P(0)=\ftilde_P(0)$ is equal to the generic rank of the cohomology sheaf $\calH^{-g}(P)$, which we will simply call the generic rank of $P$ in what follows. The functional equation therefore implies that $\deg_t \ftilde_P(t) = 2g-2$ iff the given perverse sheaf has support $\Supp(P)=A$. The examples below show that in general also $\deg_t f_P(t) = g-1$.

\medskip

The proof of theorem~\ref{thm:polynomial} also gives explicit formulae for the arising polynomials in terms of Chern-MacPherson classes~\cite{MacPhersonChern}, see theorem~\ref{thm:explicit} and \ref{thm:explicit_symmetric}. We pass to the dual abelian variety $\Adual = \Pic^0(A)$ via the Fourier transform since the latter replaces the convolution product by the intersection product~\cite{BeauvilleFourier} which is more convenient in explicit computations. The crucial input for the symmetric convolution powers is a result of Cappell, Maxim, Sch\"urmann, Shaneson and Yokura~\cite{CappellSymmetric}. One may also express  the other Schur functors via the characters of the symmetric group as follows, where $\chi_{V_\alpha}: \frakS_n \to \bbZ$ denotes the character of $V_\alpha$.

\begin{thm} \label{thm:schurintro}
Let $\alpha$ be a partition of $n$. For all clean semisimple $P\in \Pervc(A)$ with finite spectrum one has
\[
 \rtilde_P(\alpha) \;=\; \frac{1}{n!} \sum_{\sigma \in \frakS_n} \chi_{V_\alpha}(\sigma) \cdot c_P(\sigma)
\]
where the $c_P(\sigma)\in \bbZ$ are integers given by the explicit formulae in theorem~\ref{thm:schur}.
\end{thm}

Like the character values, the integers $c_P(\sigma)$ for a given $P\in \Pervc(A)$ only depend on the conjugacy class of $\sigma\in \frakS_n$. In the sequel we associate to any such class the cycle type
\[ 
 \sigma_1 \; \geq \; \sigma_2 \; \geq \; \cdots \; \geq \; \sigma_{\ell(\sigma)} \; > \; 0\]
where $\sigma_i$ denotes the length of the $i^\mathrm{th}$ cycle in a decomposition of $\sigma \in \frakS_n$ into a product of disjoint cycles. For instance we have $\ell(\sigma)=n$ iff $\sigma$ is trivial. This being said, let us consider some examples. We begin with elliptic curves:

\begin{thm} \label{thm:elliptic}
If $g=1$, then any clean semisimple $P \in \Pervc(A)$ of generic rank~$r$ and Euler characteristic $\chi = \chi(P)$ satisfies
\[
 f_P(t) \; = \; \ftilde_P(t) \; = \; r
 \quad \textnormal{\em and} \quad
 c_P(\sigma) \;=\; r \, \chi^{\ell(\sigma)-1} \, \sum_{i=1}^{\ell(\sigma)} \sigma_i^2 
 \quad \textnormal{\em for} \quad
 \sigma \in \frakS_n.
\]
\end{thm}

Note that the claim about the polynomials $f_P(t)$ and $\ftilde_P(t)$ already follows from the degree estimates in theorem~\ref{thm:polynomial} and from the above interpretation of the generic rank as the value of these polynomials at $t=0$. The formula for the numbers $c_P(\sigma)$ is less obvious but will be explained later in corollary~\ref{cor:elliptic}.

\medskip

In higher dimensions the most interesting examples arise if $P=\delta_Z \in \Perv(A)$ is the perverse intersection cohomology sheaf with support on a closed pure-dimensional subvariety~$Z\hookrightarrow A$. 
Then $P\in \Pervc(A)$ unless some irreducible component of $Z$ is invariant under translations by a positive-dimensional abelian subvariety, and we put 
\[ 
 f_Z(t)=f_P(t), 
 \quad \ftilde_Z(t)=\ftilde_P(t),
 \quad \rtilde_Z(\alpha) = \rtilde_P(\alpha)
 \quad \textnormal{and} \quad
 c_Z(\sigma) = c_P(\sigma). \]
The simplest nontrivial instance occurs for Jacobian varieties, where the arising perverse sheaves are related to the theory of special linear series on curves~\cite{WeBN}. We denote by $[h(s)]_{s^g}$ the coefficient of $s^g$ in a polynomial $h(s) \in \bbZ[s]$.

\begin{thm} \label{thm:curve}
If $A$ is the Jacobian variety of a smooth projective curve $C\hookrightarrow A$ of genus $g>1$, then 

\[
 f_C(t) \;=\; (g!  - 1) \, t^{g-1}
 \quad \textnormal{\em and} \quad
 \ftilde_C(t) \;=\; t^{g-1}
\]
and we have
\[
 c_{\hspace*{0.02em} C}(\sigma) \;=\; 
 \biggl[ \, g! \prod_{i=1}^{\ell(\sigma)} \bigl(2g-2 + \sigma_i^2 s \bigr) \, - \, \prod_{i=1}^{\ell(\sigma)} \bigl(2g-2 - s q_{\sigma_i}(s) \bigr) \, \biggr]_{s^g}
\]
where $q_{\sigma_i}(s) \in \bbZ[s]$ are polynomials given by the formula in definition~\ref{def:qpolynomial} below.
\end{thm}

The argument for Jacobian varieties generalizes directly to Prym-Tjurin varieties in the sense of~\cite{BL}, see theorem~\ref{cor:prympowers}. Another interesting case are the convolution powers of an ample divisor as in the following example~\cite{KrE6}.

 \begin{ex} \label{ex:g4}
 Let $A$ be a general principally polarized abelian variety and $\Theta \subset A$ a smooth ample divisor which defines the principal polarization. The convolution square $\delta_\Theta\conv \delta_\Theta$ defines over some open dense subset a variation of Hodge structures whose fibres are the primitive
 cohomology of the intersections $\Theta \cap (\Theta + x)\subset A$ for~$x\in A(\bbC)$. These intersections have been studied a lot in the moduli theory of abelian varieties~\cite{BeauvillePrym}~\cite{DebTorelli} ~\cite{DebarreSchottkyUpdate}~\cite{IzadiGeometricStructure}, and the higher convolutions  $\delta_\Theta \conv \cdots \conv \delta_\Theta$ provide a natural continuation of this topic. For $g = 4$, examples~\ref{ex:g4convolution} and~\ref{ex:g4symmetric} will show
\begin{align} \nonumber
 f_\Theta(t) &\;=\; 1829\, t^3 - 342\, t^2 + 58\, t,  \\ \nonumber
 \ftilde_\Theta(t) &\;=\; 52\, t^5 + 1292\, t^4 + 5049\, t^3 + 1292\, t^2 + 52\, t. 
\end{align}
In particular, on some open dense subset the Hodge modules $\delta_\Theta \circ \delta_\Theta$ and $S^2(\delta_\Theta)$ define  variations of Hodge structures of rank $58$ and $52$. Furthermore the alternating square $S^{1,1}(\delta_\Theta)$ defines a variation of Hodge structures of rank $\rtilde_\Theta(1,1)=58-52=6$ whose stalks may be identified with the lattice $E_6$ and whose monodromy group has index $\leq 2$ in $\Aut(E_6)$; this is related to the Prym morphism~\cite{KrE6}~\cite{DonagiPrym}.
\end{ex}

\section{Preliminary remarks and notations} \label{sec:notations}

The decomposition of tensor powers into Schur functors works in any $\bbQ$-linear symmetric monoidal pseudoabelian category~\cite[sect.~1.2]{DelCT}, so for any $K\in \Dbc(A)$ we may write
\[
 K^{* n} \;\cong\; \bigoplus_{\deg \alpha = n} \; S^{* \alpha}(K) \otimes_\bbQ V_\alpha
 \qquad \textnormal{with} \qquad S^{* \alpha}(K)\in \Dbc(A).
\]
If $K=P$ is a clean semisimple perverse sheaf, then the decomposition theorem and the definition of the clean perverse sheaves $P^n$ and $S^{\alpha}(P)$ from the introduction imply that
\begin{align} \nonumber
 P^{*n}  &\;\cong\; P^{n} \oplus P^{\bullet n} 
 \hspace{5.2em} \textnormal{with} \qquad P^{\bullet n} \in \rmN(A), \\ \nonumber
 S^{*\alpha}(P)  &  \;\cong\; S^{\alpha}(P) \oplus S^{\bullet \alpha}(P) 
 \hspace{2.2em} \textnormal{with} \qquad S^{\bullet \alpha}(P) \in \rmN(A).
\end{align}
So the data we are interested in becomes the difference of two parts: For $\natural \in \{*, \bullet \}$ we put
\begin{align} \nonumber
 r_P^\natural(n) &\;=\; (-1)^g \, \chi_\eta(P^{\natural n}),
 \; & Z_P^\natural(t) &\;=\; \sum\nolimits_{n\geq 0} \, r_P^\natural(n) \cdot t^n \\ \nonumber
 \rtilde_P^\natural(\alpha) &\;=\; (-1)^g \, \chi_\eta(S^{\natural \alpha}(P)),
 \; & \Ztilde_P^\natural(t) &\;=\; \sum\nolimits_{n\geq 0} \, \rtilde_P^\natural(n) \cdot t^n 
\end{align}
where 
$
 \chi_\eta(-) 
$
denotes the Euler characteristic of the stalk cohomology at the generic point of the abelian variety. Then
\[
 Z_P \;=\; Z_P^*  \;-\;  Z_P^\bullet, 
 \qquad
 \Ztilde_P \;=\; \Ztilde_P^* \;-\; \Ztilde_P^\bullet
 \qquad \textnormal{and} \qquad
 \rtilde_P \;=\; \rtilde_P^* \;-\; \rtilde_P^\bullet,
\]
and the results from the introduction will be shown separately for each term.

\medskip

To do this in a uniform way, let $\rmD(A) \subseteq \Dbc(A)$ be the full additive subcategory of all sheaf complexes that decompose as a direct sum $\bigoplus_{n\in \bbZ} P_n[n]$ for semisimple perverse sheaves $P_n \cong P_{-n}$ with finite spectrum. From the decomposition theorem and the relative hard Lefschetz theorem one sees that this subcategory is stable under the convolution product. Consider the associated Grothendieck ring
\[\rmK(A) \;=\; \rmK_0(\rmD(A), \oplus, \, * \, )  \]
where the underlying Grothendieck group is taken only with respect to direct sums, not distinguished triangles: We keep the information on the degrees in which a complex sits. To compute the above terms with the superscript $\natural = \bullet$ resp.~$\natural = *$ we consider ring homomorphisms 
\[
\xymatrix@M=0.5em@R=1.8em@C=3.5em{ 
  &&& \rmK(A) \ar[dl]_-{h_s} \ar[dr]^-{\gamma_s} &&& \\
  && \ar@{}[ull]|{\textnormal{(for $\natural = \bullet$)}}\bbQ[s]  &  & \rmH^{2\bullet}(\Adual, \bbQ)[s]  \ar@{}[urr]|{\textnormal{(for $\natural = *$)}} && 
}
\]
that will be defined, respectively, in section~\ref{sec:betti} via signed Poincar\'e polynomials and in section~\ref{sec:chern} via Chern-MacPherson classes. Here $\rmH^{2\bullet}(\Adual, \bbQ) = \bigoplus_{n\in \bbN_0} \rmH^{2n}(\Adual, \bbQ)$ is the even degree cohomology ring. For power series $f(s) = \sum_{n\geq 0} r_n s^n \in R[[s]]$ we write
\[
 \bigl[ f(s) \bigr]_{s^n} \;=\; 
 \begin{cases}
 r_n & \textnormal{if} \; R = \bbQ, \\
 ev(r_n) & \textnormal{if} \; R = \rmH^{2\bullet}(\Adual, \bbQ),
 \end{cases}
\]
where $ev: \rmH^{2\bullet}(\Adual, \bbQ) \rightarrow \rmH^{2g}(\Adual, \bbQ) \rightarrow \bbQ$ is the evaluation on the fundamental class.

\section{Signed Poincar\'e polynomials and negligible complexes} \label{sec:betti}

To control the negligible direct summands which occur in convolution powers of a complex $K\in \rmD(A)$, we consider the Betti numbers $h^n(K) =\dim_\bbC \rmH^n(A, K)$ and the signed Poincar\'e polynomial

\[ 
 b_x(K) \;=\; \sum_{n\in \bbZ} (-1)^n h^n(K) \, x^n \;\in \; \bbZ[x^{\pm 1}]. \] 
The signs are required for the compatibility with the pre-$\lambda$-structures in lemma~\ref{lem:lambdahom} below.
We also consider the twists $K_\varphi = K\otimes_\bbC L_\varphi$ by rank one local systems~$L_\varphi$ given by characters
\[ \varphi \; \in \; \Pi(A) \;=\; \Hom_\bbZ(\pi_1(A, 0), \bbC^*). \]
For a semisimple perverse sheaf $P\in \Perv(A)$, the spectrum $\scrS(P)$ is finite iff $b_x(P_\varphi)$ is a constant polynomial for all but finitely many~$\varphi$, indeed
\[ 
 \scrS(P) \;=\; \bigl\{\varphi \in \Pi(A) \mid \exists n\neq 0: \; \rmH^n(A, P_\varphi)\neq 0 \bigr\}. 
\]
We define the spectrum $\scrS(K)$ of a complex $K\in \rmD(A)$ to be the union of the spectra of its perverse cohomology sheaves. In what follows it will be convenient to replace the variable $x$ in the signed Poincar\'e polynomial by a new variable: The hard Lefschetz theorem and the definition of~$\rmD(A)$ imply that $h^n(K)=h^{-n}(K)$ for all $n\in \bbZ$, so we get a factorization
\[
\xymatrix@M=0.5em@C=3em{
 \rmK(A) \ar@{..>}[dr]_-{\exists! h_s} \ar[rr]^-{b_x} && \bbZ[x^{\pm 1}] \\
 & \bbZ[s] \ar@{^{(}->}[ur]_-\iota & 
}
\]
with $\iota$ given by
\[
 \iota(s) \;=\;  2 - x^{-1} - x.\medskip
\]
Note that the value of the polynomial $h_s(K)$ at $s=0$ coincides with the Euler characteristic $\chi(K)$. Furthermore our new variable has the property $s^g = h_s(\delta_A)$, which results in the following formula.

\begin{prop} \label{prop:negligible}
If $K=P\oplus N$ where $P\in \Pervc(A)$ is a clean semisimple perverse sheaf and $N\in \rmD(A)$ is a negligible complex with finite spectrum, then $h_s(K_\varphi)$ is a polynomial of degree
\medskip
$$
 \deg_s h_s(K_\varphi) \;<\;g \quad \textnormal{\em for all} \quad \varphi \; \in \; \Pi(A) \setminus \scrS(N),
$$
and \medskip
\[
 (-1)^g \, \chi_\eta(N) \;\;\;=\; \sum_{\varphi \in \scrS(N)} \bigl[ h_s(K_\varphi) \bigr]_{s^g}
 \;=\; \sum_{\varphi \in \Pi(A)} \bigl[ h_s(K_\varphi) \bigr]_{s^g}.
\]
\end{prop}

{\em Proof.} We first claim $\rmH^n(A, P)=0$ for $|n|\geq g$. Indeed, for $|n|>\dim \Supp(P)$ this vanishing holds by the standard estimate for the cohomological dimension of direct image functors with respect to the perverse $t$-structure~\cite[cor.~6.0.5]{SchuermannTopology}, so we may assume that $|n|=\dim \Supp(P) =g$. We may also assume that the perverse sheaf $P$ is simple and hence arises as the intermediate extension $P = j_{!*}( L[g])$ of a simple local system $L$ on some open dense subset $j: U\hookrightarrow A$. The excision sequence for the closed complement $i: Z = A\setminus U \hookrightarrow A$ then takes the form
$$
 0 \;=\; \rmH^{-g}(Z, i^!(P)) \;\to\; \rmH^{-g}(A, P) \;\to\; \rmH^{-g}(Z,j^*(P)) \;=\; \rmH^0(U, L) \;\to\; \cdots
$$
where the zero on the left is again due to the standard estimate for the cohomological dimension since $i^!(P)$ is semiperverse and $\dim Z < g$. But we also have $\rmH^0(U, L) = 0$ since otherwise 
the simple local system $L$ would have to be trivial, contradicting our assumption that $P$ is clean. Altogether this shows $\rmH^{-g}(A, P)=0$, and then also  $\rmH^{g}(A, P)=0$ by the same argument applied to the Verdier dual.

\medskip

 It follows that for any clean semisimple perverse sheaf $P\in \Pervc(A)$ the signed Poincar\'e polynomial has degree $\deg_s h_s(P) < g$, so we may assume from the start that $K=N$ is negligible. Any simple negligible perverse sheaf with finite spectrum has the form $L_\varphi[g]$ for some $\varphi \in \Pi(A)$ by~\cite[prop.~10.1]{KrWVanishing}, and $\rmH^\bullet(A, L_\varphi)=0$ for all $\varphi \neq 0$. Since the claim of the proposition is additive under direct sums and invariant under character twists, we may thus assume $N =  V^\bullet \otimes_\bbC \delta_A$ for a bounded complex $V^\bullet$ of vector spaces. Then
\[
 b_x(N) \;=\; b(x)\cdot b_x(\delta_A) \quad \textnormal{with} \quad
 b(x) \;=\; \sum\nolimits_{n\in \bbZ} (-1)^n \, \dim_\bbC(V^n) \cdot x^n 
\]
and our claim follows by direct inspection: Writing $b(x)=\iota(h(s))$ with $h(s)\in \bbZ[s]$, we have $(-1)^g \chi_\eta(N) = b(1) = h(0)$. \qed

\medskip

For the symmetric powers, recall that a {\em pre-$\lambda$-structure} on a commutative ring~$R$ with unit~$1$ is by definition a  homomorphism
$$
 S_t: \quad  (R, \, + \, )  \; \longrightarrow \; (1+tR[[t]], \, \cdot \, )
$$
from the additive group of the ring to the multiplicative group of monic power series such that $S_t(a) \equiv 1 + a \, t \;\, \textnormal{mod}\;\,  t^2 R[[t]]$ for all $a\in R$, see~\cite{KnutsonLambda}.

\begin{ex}
On the Grothendieck ring $\rmK(A)$ the symmetric convolution powers define a pre-$\lambda$-structure via
\[
 S_t(K) \;=\; 1 \; + \; \sum_{n=1}^\infty \; S^{*n}(K) \cdot t^n \quad \textnormal{for} \quad K\in \rmK(A).
\]
The same construction works for the Grothendieck ring of any  $\bbQ$-linear symmetric monoidal pseudoabelian category, and we call it the symmetric pre-$\lambda$-structure.
\end{ex}

\begin{ex}\label{ex:betti-lambdastructure}
We will identify $\bbZ[x^{\pm 1}]$ with the Grothendieck ring of the category of~$\bbZ$-graded finite-dimensional complex vector spaces $V$ by the assignment

\[ 
 V \;=\; \bigoplus_{n\in \bbZ} \; V_n \; \mapsto \; \sum_{n\in \bbZ} \, (-1)^n \dim_\bbC(V_n) \, x^n \; \in \; \bbZ[x^{\pm 1}]. 
\] 
With the usual sign rule for graded tensor products, the symmetric pre-$\lambda$-structure then becomes
\[
 S_t \Bigl(\, \sum\nolimits_{n\in \bbZ} a_n x^n \Bigr) \;=\; \prod_{n\in \bbZ} \bigl( 1 - x^n t \bigr)^{-a_n}
\]
where on the right each factor has to be expanded as a power series.
\end{ex}

\begin{lem} \label{lem:lambdahom}
On the Grothendieck ring the signed Poincar\'e polynomial induces a ring homomorphism
 $$b_x: \quad \rmK(A) \; \longrightarrow \; \bbZ[x^{\pm 1}]$$ 
which is compatible with the pre-$\lambda$-structures in the above two examples.
\end{lem}

{\em Proof.} The K\"unneth formula implies that the signed Poincar\'e polynomial induces a ring homomorphism from $\rmK(A)$ to $\bbZ[x^{\pm 1}]$. More precisely, the definition of the commutativity constraint for the convolution product implies that $\rmH^\bullet(A, -)$ is a tensor functor from the symmetric monoidal category $(\rmD(A), *)$ to the Tannakian category of finite-dimensional $\bbZ$-graded complex vector spaces~\cite[lemma 8]{WeBN}. Hence the compatibility of $b_x$ with the above pre-$\lambda$-structures follows from our choice of signs in the definition of $b_x$.
  \qed

\medskip

\begin{dfn} \label{def:qpolynomial}
We define $q_n(s)\in \bbZ[s]$ by
$ x^n + x^{-n} - 2 = \iota(s\, q_n(s))$.
Explicitly one may check that
\[
 q_n(s) \;=\; - \sum_{e=1}^n 
 \tbinom{n + e}{n - e} \cdot \tfrac{2n}{n+e}
  \cdot (-s)^{e-1}
\]
but this will not be relevant now. We put $\nu_n(K) = (-1)^{n+1} h^n(K)$ for $K\in \rmD(A)$.
\end{dfn}

\begin{cor} \label{cor:bettiformula}
For any $K\in \rmD(A)$,
\[
 \sum_{n=0}^\infty h_s(S^{*n}(K)) \cdot t^n \;=\; 
 (1-t)^{-\chi(K)} \cdot \prod_{n>0} \Bigl(1 - \frac{q_n(s) \cdot st}{(1-t)^2} \, \Bigr)^{\nu_n(K)}
\]
\end{cor}

{\em Proof.} By lemma~\ref{lem:lambdahom} the series on the left is the preimage of~$S_t(b_x(K))$ under the variable transform~$\iota: \bbZ[s] \hookrightarrow \bbZ[x^{\pm 1}]$, so this is really only a statement about the pre-$\lambda$-structure which is induced on the subring $\bbZ[s]$ via this transform. Both sides of the claim are multiplicative with respect to sums of polynomials, so it suffices to show  $S_t(x^n + x^{-n}) = (1-t)^{-2} \cdot \iota(1 - \tfrac{q_n(s) \cdot st}{(1-t)^2})^{-1}$ for all $n\in \bbN$. But this follows via a direct computation from the formula for the pre-$\lambda$-structure in example~\ref{ex:betti-lambdastructure}. \qed

\section{Chern-MacPherson classes and generic stalks} \label{sec:chern}

Having dealt with the negligible terms, we now determine the Euler characteristic of the generic stalk cohomology of a convolution product. We denote by $\rmF(A)$ the group of constructible functions $A(\bbC) \rightarrow \bbZ$, such as for instance the stalkwise Euler characteristic

\[
  \chiloc(K): \quad A(\bbC) \; \longrightarrow \;\bbZ, 
  \quad s \; \mapsto \;  \sum_{n\in \bbZ} \, (-1)^n \dim \calH^n(K)_s
\]
of a constructible complex $K\in \Dbc(A)$. This Euler characteristic only depends on the class of the complex in the Grothendieck ring $\rmK_0(A) = \rmK_0(\Dbc(A), \oplus, *)$. In this section we will always work on the level of this Grothendieck ring, but the ring homomorphism
$$
 \rmK(A) \;=\; \rmK_0(\rmD(A), \oplus, *) \;\longrightarrow \; \rmK_0(A) \;=\; \rmK_0(\Dbc(A), \oplus, *)
$$
induced by the inclusion functor $\rmD(A) \subset \Dbc(A)$ will allow to read our results also in the Grothendieck ring from the previous sections. Let $\rmH^{2\bullet}(-)=\rmH^{2\bullet}(-, \bbQ)$ denote the even degree cohomology ring with rational coefficients, and similarly for homology. Consider the composite map
\[
 \gamma: \quad
 \xymatrix@M=0.5em@C=3em{
 \rmK_0(A) \ar[r]^-{\chiloc} 
 & \rmF(A) \ar[r]^-c
 & \rmH_{2\bullet}(A) \; \cong \; \rmH^{2g-2\bullet}(A) \ar[r]^-F
 & \rmH^{2\bullet}(\Adual) 
}
\]
where $c$ denotes the Chern-MacPherson class as defined in~\cite{MacPhersonChern} and where $F$ is the Fourier transform~\cite{BeauvilleFourier}; the identification of homology and cohomology in the middle is the Poincar\'e duality isomorphism and will be suppressed in what follows.

\begin{lem}
The above map is a ring homomorphism 
$\gamma: \rmK_0(A) \longrightarrow \rmH^{2\bullet}(\Adual)$ 
with respect to the convolution product $*$ on the Grothendieck ring and the intersection product $\, \cdot \,$ on the cohomology ring.
\end{lem}

{\em Proof.} The additivity of $\gamma$ is clear from the definitions, so we only need to check multiplicativity. Recall~\cite{BeauvilleFourier} that the convolution product $\rmH^{\bullet}(A) \times \rmH^{\bullet}(A) \to \rmH^{\bullet}(A)$, defined by $(\alpha, \beta) \mapsto \alpha * \beta = a_*(\alpha \boxtimes \beta)$, corresponds under the Fourier transform to the intersection product on the dual abelian variety: $F(\alpha * \beta) = F(\alpha) \cdot F(\beta)$. This reduces our claim to the commutativity of the following diagram.
\[
\xymatrix@R=3.5em@C=3.2em@M=0.5em{
\rmK_0(A) \times \rmK_0(A) \ar[r]^-{\chiloc} \ar[d]^{\boxtimes} 
& \rmF(A) \times \rmF(A) \ar[r]^-{c} \ar[d]^{\boxtimes} 
& \rmH^{2\bullet}(A) \times \rmH^{2\bullet}(A) \ar[d]^{\boxtimes} 
\\
\rmK_0(A\times A) \ar[r]^-\chiloc \ar[d]^{Ra_*} 
& \rmF(A\times A) \ar[r]^-{c} \ar[d]^{a_*} 
& \rmH^{2\bullet}(A\times A) \ar[d]^{a_*} 
 \\
\rmK_0(A) \ar[r]^-\chiloc 
& \rmF(A) \ar[r]^-{c} 
& \rmH^{2\bullet}(A)
}
\]
\noindent
Here the pushforward of constructible functions is defined by the fibrewise weighted Euler characteristic~\cite{MacPhersonChern}. Note that this pushforward is compatible with the proper direct image of constructible sheaf complexes~\cite[sect.~2.3]{SchuermannTopology}. The compatibility of $c$ with proper direct images follows from the naturality of Chern-MacPherson classes and the compatibility with $\boxtimes$ has been established in~\cite{Kwiecinski_Formule} and~\cite{KY_Product}.
\qed

\medskip

For $K\in \rmK_0(A)$ we denote by $\gamma_n(K) \in \rmH^{2n}(\Adual)$ the degree $2n$ part of $\gamma(K)$ and consider the polynomial
\[
 \gamma_s(K) \;=\; \sum_{n=0}^g \; \gamma_n(K) \cdot s^n \;\in\; \rmH^{2\bullet}(\Adual)[s].
\]

\begin{lem} \label{lem:euler}
For all $K\in \rmK_0(A)$ one has $(-1)^g \, \chi_\eta(K) = \gamma_g(K) = [ \gamma_s(K) ]_{s^g}$.
\end{lem}

{\em Proof.} This amounts to the commutativity of the following diagram, where $ev_\eta$ is the evaluation of constructible functions at the generic point and $pr$ denotes the natural projections.
\[
\xymatrix@R=2.5em@C=3em@M=0.5em{
& \rmF(A) \ar[r]^-{c} \ar[dd]^-{ev_\eta} 
& \rmH^{2\bullet}(A) \ar[r]^-F \ar[d]^-{pr} 
& \rmH^{2\bullet}(\hat{A}) \ar[d]^-{pr}
\\
\rmK_0(A) \ar[ur]^-\chiloc \ar[dr]_-{\chi_\eta}
& 
& \rmH^0(A) \ar@{=}[d] \ar[r]^-{F} 
& \rmH^{2g}(\hat{A})  \ar[d]^-{ev} 
\\
& \bbZ \ar@{^{(}->}[r] 
& \bbQ \ar[r]^-{(-1)^g} 
& \bbQ
}
\]
For the sign $(-1)^g$ in the bottom right part we refer to~\cite[prop.~1]{BeauvilleFourier}. \qed

\medskip

Similarly one sees that $\gamma_0(K) = \chi(K)$ by functoriality of the Chern-MacPherson class, since $A$ is connected.
There is also a simple generating series for the symmetric convolution powers in this context; however, here it will be easier to work not with pre-$\lambda$-structures but instead with the corresponding Adams operations which give an expression as an exponential series. Recall that 
\[
 \sum_{r=1}^\infty \; r^m t^r \;=\; \frac{t \, p_m(t)}{(1-t)^{m+1}}
 \quad \textnormal{for all} \quad m \; > \; 0,
\]
where the $p_m(t) \in \bbZ[t]$ are monic polynomials of degree $\deg_t p_m(t) = m-1$ which are also known as the {\em Eulerian polynomials}~\cite{HirzebruchEulerian}.

\begin{prop} \label{prop:symmetric}
For any $K\in \Dbc(A)$, 
\[
\sum_{n=0}^\infty \gamma (S^{*n}(K)) \cdot t^n 
\;=\; 
(1-t)^{-\chi(K)} \cdot 
\prod_{i=1}^g 
\exp \Bigl(  \gamma_i(K) \cdot \frac{tp_{2i-1}(t)}{(1-t)^{2i}}\,
\Bigr).
\]
\end{prop}

{\em Proof.} For $n\in \bbN$, let us denote by $q_n: A^n \twoheadrightarrow A^{(n)} = A^n/\frakS_n$ the projection onto the $n$-fold symmetric product of the abelian variety. Then the symmetric group $\frakS_n$ naturally acts on $q_{n*}(K \boxtimes \cdots \boxtimes K)$, and since $\Dbc(A^{(n)})$ is a $\bbQ$-linear pseudoabelian category, we may define the symmetric product
$$
 K^{(n)} \;=\; (q_{n*}(K\boxtimes \cdots \boxtimes K))^{\frakS_n} \;\in\; \Dbc(A^{(n)})
$$
by applying the Schur functor that projects onto the invariants under the symmetric group. The functoriality of this projector implies 
$$ S^{*n}(K) \;=\; a_{n*}(K^{(n)}) \quad \textnormal{for the addition morphism} \quad a_n: \; A^{(n)} \;\longrightarrow\; A. $$
We now use a generating series for the Chern-MacPherson classes of the symmetric products $K^{(n)}$ that has been obtained in~\cite[th.~1.7]{CappellSymmetric} via the Adams operations. The generating series in loc.~cit.~is expressed via the pushforwards $d_{r*}(K)$ for $r\in \bbN$, where $d_r: A \hookrightarrow A^{(r)}$ denotes the diagonal embedding. To translate this result to our present setting, note that the composite $a_r\circ d_r = [r]$ is the multiplication by $r$ endomorphism on the abelian variety $A$. The Fourier transform converts the pushforward $[r]_*$ on $H_{2\bullet}(A)$ into the pull-back $[r]^*$ on $H^{2\bullet}(\Adual)$ by~\cite[prop.~3(iii)]{BeauvilleFourier}, so we obtain the generating series                                                                                                                                                                                                                   %
\begin{equation} \label{eq:exponential}
 \sum_{n=0}^\infty \gamma(S^{*n}(K)) \cdot t^n 
 \;=\; 
 \exp \Bigl( \, \sum_{r = 1}^\infty \; [r]^*(\gamma(K)) \cdot \frac{t^r}{r} \, \Bigr) .
\end{equation}
Now $[r]^*$ acts as multiplication by $r^m$ on $H^m(\Adual)$, so our claim follows by inserting the value~$\gamma_0(K)=\chi(K)$ and by using that $\exp ( \chi(K) \sum_{r=1}^\infty t^r/r ) = (1-t)^{-\chi(K)}$.
\qed

\medskip

\section{Generating series for convolution powers} \label{sec:series}

Recalling the decomposition $Z_P(t)=Z_P^*(t) - Z_P^\bullet(t)$ from section~\ref{sec:notations} and putting together the information on signed Poincar\'e polynomials and Chern-MacPherson classes from above, we may now prove the first part of theorem~\ref{thm:polynomial} as follows.

\begin{thm} \label{thm:explicit}
Let $P\in \Pervc(A)$ be a clean semisimple perverse sheaf and $\chi = \chi(P)$ its Euler characteristic. Let $\natural \in \{ *, \bullet\}$, and for $\natural = \bullet$ assume furthermore that the spectrum $\scrS = \scrS(P)$ is finite. Then

\[
 Z_P^\natural(t) \;=\; \frac{t f_P^\natural(t)}{(1-\chi t)^{g+1}}
\]
where $f_P^\natural(t) \in \bbZ[t]$ is given by
$$
f_P^\natural(t) \;\; = \;\;
\begin{dcases}
\quad \sum_{n=1}^g  \; \bigl[ (\gamma_s(P)-\chi)^n  \bigr]_{s^g} \! \cdot t^{n-1}(1-\chi t)^{g-n}  & \textnormal{\em if} \quad \natural \;=\; *, \\
\quad \sum_{n=1}^g \sum_{\varphi \in \scrS} \; \bigl[ (h_s(P_\varphi)-\chi)^n \bigr]_{s^g} \! \cdot t^{n-1}(1-\chi t)^{g-n}  & \textnormal{\em if} \quad \natural \;=\; \bullet.
\end{dcases}
$$
\end{thm}

{\em Proof.} We have $\gamma_s(P^{*n})=(\gamma_s(P))^n$ and $h_s(P^{*n})=(h_s(P))^n$ since $\gamma_s$ and~$h_s$ are ring homomorphisms. The same also holds after character twists since twisting by a character is a tensor functor~\cite[prop.~4.1]{KrWVanishing}. Hence proposition~\ref{prop:negligible} and lemma~\ref{lem:euler} show 
\[ 
 r_P^*(n) \;=\; \bigl[ (\gamma_s(P))^n \bigr]_{s^g} 
 \quad \textnormal{and} \quad 
 r_P^\bullet(n) \;=\; \sum\nolimits_{\varphi \in \scrS} \bigl[ (h_s(P_\varphi))^n \bigr]_{s^g}. 
\]
The result now follows from a simple remark about power series: Let $h(s) \in R[[s]]$ be a power series over a commutative ring $R$ and $\chi =h(0)$. Then we claim that for any $g\in \bbN$ the identity
\[
 (1-\chi t)^{g+1} \, \sum_{n=1}^\infty \bigl[ (h(s))^n \bigr]_{s^g}  t^n \;=\; 
 \sum_{n=1}^g \, \bigl[(h(s) - \chi)^n\bigr]_{s^g} \cdot t^n \, (1-\chi t)^{g-n} 
\]
holds. Indeed this follows from
\[
 (1-\chi t) \, \sum_{n=0}^\infty \, (h(s) t )^n 
 \;=\; \frac{1-\chi t}{ \; 1 - h(s) t \; } 
 \;=\; \sum_{n=0}^\infty \, \Bigl( \frac{h(s) - \chi}{1-\chi t} \Bigr)^n \cdot t^n  
\]
because the coefficient of $s^g$ in $(h(s)-\chi)^n$ vanishes for $n>g$.
\qed

\medskip

\begin{ex} \label{ex:g4convolution}
Let $i: \Theta \hookrightarrow A$ be a smooth ample divisor whose class $\theta \in \rmH^2(A)$ defines a principal polarization. The smoothness implies that $\delta_\Theta = i_*(\bbC_\Theta)[g-1]$ and hence $c(\delta_\Theta) = (-1)^{g-1} i_*(c(\bbC_\Theta))$. By definition the Chern-MacPherson class of the constant function on a smooth variety is Poincar\'e dual to the total Chern class~$c^*$ of the variety, so the class $i_*(c(\bbC_\Theta))$ in $\rmH_{2\bullet}(A)$ is the Poincar\'e dual of $i_!(c^*(\Theta))$ in $\rmH^{2\bullet}(A)$. Since the tangent bundle to abelian varieties is trivial, this latter class may be computed from the adjunction formula to be $i_!(c^*(\Theta)) = \theta\cdot (1+\theta)^{-1}$ and altogether we get
$$ c(\delta_\Theta) \;=\; \sum_{k=1}^g \; (-1)^{k+g}\, \theta^k \quad\textnormal{in}\quad \rmH^{2\bullet}(A).
$$
Now recall that we assumed $\theta$ defines a principal polarization. If $\tdual \in \rmH^2(\Adual)$ denotes the class of the dual principal polarization, we get from the formulae for the Fourier transform in~\cite[lemme~1]{BeauvilleFourier} that
\[
 \gamma(\delta_\Theta) \;=\;  \sum_{k=1}^{g} \, \tfrac{k!}{(g-k)!} \cdot \tdual^{g-k} \; \in \; \rmH^{2\bullet}(\Adual).
\]
On the other hand the weak Lefschetz theorem together with $\chi(\delta_\Theta)=g!$ determines the cohomology of the theta divisor, and one may check that this leads to the signed Poincar\'e polynomial
\[
 h_s(\delta_\Theta) \;=\; g! \;+\; \sum_{k=1}^{g-1} \, \tbinom{2k}{k} \cdot \tfrac{1}{k+1} \cdot s^{g-k}\;\in\; \bbZ[s].
\]
For $g=4$ this implies

\[
 f_\Theta^*(t) \;=\; 4! \, (432\, t^3 - 36\, t^2 + 3\, t ), 
 \qquad
 f_\Theta^\bullet(t) \;=\; 8539\, t^3 - 522\, t^2 + 14\, t.
\]
\end{ex}

\medskip

For the symmetric powers we again consider the  Eulerian polynomials $p_m(s)$ of section~\ref{sec:chern} and the polynomials $q_n(s)$ from definition~\ref{def:qpolynomial}.

\begin{thm} \label{thm:explicit_symmetric}
Let $P\in \Pervc(A)$ be a clean semisimple perverse sheaf and $\chi=\chi(P)$ its Euler characteristic. Let $\natural \in \{*, \bullet\}$, and for $\natural = \bullet$ assume furthermore that the spectrum $\scrS = \scrS(P)$ is finite. Then

\[
 \Ztilde_P^\natural(t) \;=\; \frac{t\ftilde_P^\natural(t)}{(1-t)^{2g+\chi}} 
\]
where $\ftilde_P(t) \in \bbZ[t]$ is an integral polynomial of degree at most $2g-2$ satisfying the functional equation
\[
 \ftilde_P^\natural(t) \;=\; t^{2g-2}\, \ftilde_P^\natural(1/t).
\]
Explicitly,
$$
 \ftilde_P^\natural(t) \;\; = \;\;
 \begin{dcases}
 \;\;  t^{-1} \Bigl[ \, \prod_{i = 1}^g  \exp \bigl(  \gamma_i(P)\,t p_{2i-1}(t)\,  s^i\,
  \bigr) \, \Bigr]_{s^g} 
  & \textnormal{\em if} \quad \natural \;=\; *,
  \\
 \;\;  \sum_{\varphi \in \scrS} t^{-1} \Bigl[ \, \prod_{n>0}  \bigl( 1-q_n((1-t)^2 \, s)\, st \, \bigr)^{\nu_n(P_\varphi)} \, \Bigr]_{s^g}
 & \textnormal{\em if} \quad \natural \;=\; \bullet,
 \end{dcases}
$$
where in the last expression the exponents are defined by $\nu_n(P_\varphi)=(-1)^{n+1}h^n(P_\varphi)$.
\end{thm}

{\em Proof.} Take the coefficient of $s^g$ in the series in corollary~\ref{cor:bettiformula} and the projection of proposition~\ref{prop:symmetric} onto the top cohomology. In both cases a factor $t/(1-t)^{2g}$ may be extracted from the product on the right hand side, and the remaining expression is precisely the polynomial $\ftilde_P^\natural(t)\in \bbZ[t]$ in the explicit form claimed above. It remains to establish the functional equation for these explicit polynomials. For $\natural = *$ this follows from the fact that the coefficients of the Eulerian polynomials satisfy the recursion in~\cite[eq.~(7.b)]{HirzebruchEulerian}, leading to a generalized Pascal triangle whose symmetry amounts to the required property $p_{2i-1}(t) = t^{2i} p_{2i-1}(1/t)$. For $\natural = \bullet$ the functional equation follows from the identity
\medskip
$$ ((1-t)^2 s)^i (st)^k \;=\; t^{2(i+k)}((1-1/t)^2 s)^i (s/t)^k 
\medskip
$$
for $i, k\in \bbN$, applied to the given explicit formula for the polynomial $f_P^\bullet(t)$.
\qed

\medskip

\begin{ex} \label{ex:g4symmetric}
Consider again a smooth ample divisor $\Theta \subset A$ whose class defines a principal polarization. For $g=4$, inserting 
\begin{align} \nonumber
p_1(t) &\;=\; 1 
& q_1(s) &\;=\; -1 \\ \nonumber
p_3(t) &\;=\; t^2 + 4t + 1  
& q_2(s) &\;=\; s-4 \\ \nonumber
p_5(t) &\;=\; t^4 + 26t^3 + 66t^2+26t +1 
& q_3(s) &\;=\; -s^2 +6s-9 
\end{align}
and the Chern-MacPherson classes and signed Poincar\'e polynomials of example~\ref{ex:g4convolution} in the theorem, we get
\begin{align} \nonumber
 \ftilde_\Theta^*(t) &\;=\; 36\, t^5 + 1152\, t^4 + 4824\, t^3 + 1152\, t^2 + 36\, t \\ \nonumber
 \ftilde_\Theta^\bullet(t) &\;=\; -16\, t^5 - 140\, t^4 - 225\, t^3 - 140\, t^2 - 16\, t
\end{align}
which in particular verifies the claim for the symmetric powers in example~\ref{ex:g4}.
\end{ex}

\section{Other Schur functors}

The above methods can also be used to compute the generic ranks $\rtilde_P(\alpha)$ as follows. Let $\scrC$ be a $\bbQ$-linear pseudoabelian symmetric monoidal category in the sense of~\cite[sect.~1.2]{DelCT}. 
For any object $P\in \scrC$ and $n\in \bbN$ we have an $\frakS_n$-equivariant decomposition
$$
 P^{\otimes n} \;\;\cong\; \bigoplus_{\deg \alpha = n} \; S^\alpha(P) \otimes_\bbQ V_\alpha
$$
where $\alpha$ runs through all partitions of degree $n$ and where $S^\alpha: \scrC \rightarrow \scrC$ are the Schur functors corresponding to the irreducible representations $V_\alpha \in \Rep_\bbQ(\frakS_n)$ of the symmetric group. The characters of all these representations are integer valued functions $\chi_{V_\alpha}: \frakS_n \to \bbZ$, so for each $\sigma\in \frakS_n$ we obtain a trace map~\cite[sect.~3]{MaximSchuermannTwisted}

\[
tr(\sigma): \quad
 \rmK_0 \; \longrightarrow \; \rmK_0,
 \quad
 [P] \; \mapsto \! \sum_{\deg \alpha = n} \chi_{V_\alpha}(\sigma) \cdot [S^\alpha(P)].
\]
on the Grothendieck ring $\rmK_0 = \rmK_0(\scrC, \oplus, \otimes)$.

\begin{ex} \label{ex:adams}
For the $n$-cycle $z_n = (12\cdots n) \in \frakS_n$ one sees by~\cite[I.4]{KnutsonLambda} \cite[cor.~1.8]{AtiyahPower} that the trace map
\[ tr(z_n) \;=\; \Psi^n \] 
is equal to the $n$-th Adams operation for the symmetric pre-$\lambda$-structure on $\rmK_0$. The Adams operations satisfy 
\[
 S_t(P) \;=\; \exp \Bigl(\, \sum\nolimits_{n > 0} \Psi^n(P) \, \frac{t^n}{n} \,\Bigr).
\]
Furthermore, in our case each $\Psi^n$ is a ring homomorphism since by~\cite[lemma~4.1]{Heinloth} the symmetric pre-$\lambda$-structure is in fact a $\lambda$-structure, i.e.~special in the terminology of loc.~cit. We require the following two cases: \vspace*{0.15em}
\begin{enumerate}
\item The Grothendieck ring $\rmK_0 = \bbZ[x^{\pm 1}]$ of the category of graded vector spaces in example~\ref{ex:betti-lambdastructure} has Adams operations given by
$\Psi^n(x^i) = x^{in}$ for $i\in \bbZ$. \smallskip
 
\item For $\rmK_0 = \rmK(A)$ the formula~\eqref{eq:exponential} shows that the Adams operations and the Fourier transformed Chern-MacPherson classes satisfy $\gamma \circ \Psi^n = [n]^* \circ \gamma$.
\end{enumerate}
\end{ex}
\noindent
As in the introduction we attach to $\sigma \in \frakS_n$ the partition $(\sigma_1, \sigma_2, \dots, \sigma_{\ell(\sigma)})$ of~$n$, where $\sigma_1 \geq \sigma_2 \geq \cdots \geq \sigma_{\ell(\sigma)} > 0$ are the lengths of the cycles in a disjoint cycle decomposition of the given permutation. With these notations we have the following abstract version of the averaging property in \cite[sect.~3]{CappellSymmetric}.

\begin{lem} \label{lem:projector}
The image of any object $P\in \scrC$ under the Schur functor $S^\alpha: \scrC \to \scrC$ has in $\rmK_0$ the class 
\[
 [S^\alpha(P)] \;=\; 
 \frac{1}{n!} \sum_{\sigma \in \frakS_n} \chi_{V_\alpha}(\sigma) 
 \cdot \prod_{r=1}^{\ell(\sigma)} 
 \Bigl(
 \Psi^{\sigma_r}[P]
 \Bigr).
\]
\end{lem}

{\em Proof.} If $a,b\in \bbN$, then for any $g\in \frakS_{a}$, $h\in \frakS_{b}$ with product $g\times h\in \frakS_{a+b}$ we have a commutative diagram involving the diagonal map $\Delta$:

\[
\xymatrix@R=3em@C=5em@M=0.6em{
 \rmK_0 \ar[r]^-{tr(g\times h)} \ar[d]_-{\Delta} 
& \rmK_0
\\
\rmK_0 \otimes_\bbZ \rmK_0 \ar[r]^-{tr(g)\otimes tr(h)}
& \rmK_0 \otimes_\bbZ \rmK_0  \ar[u]_-{\textnormal{``$\, \cdot\, $''}} 
}
\]
This follows from the Littlewood-Richardson rule for the decomposition of tensor products~\cite[prop.~1.6]{DelCT}. For the disjoint cycle decomposition of $\sigma \in \frakS_n$, example~\ref{ex:adams} then implies 
\[ 
 tr(\sigma)[P] \;=\; \prod_{r=1}^{\ell(\sigma)} \Bigl( \Psi^{\sigma_r}[P] \Bigr).
\]
Hence it only remains to note that the endomorphism
$ 
 \frac{1}{n!} \sum_{\sigma \in \frakS_n} \chi_{V_\alpha}(\sigma) \cdot tr(\sigma)(-)
$
is given by $[P] \mapsto [S^\alpha(P)]$. This follows from the orthonormality relations for the characters of the symmetric group. \qed

\medskip

As an immediate corollary we obtain the formula in theorem~\ref{thm:schurintro}. For $n\in \bbN$ we denote by
\[ 
 [n]^*: \; \rmH^\bullet(\Adual) \longrightarrow \rmH^\bullet(\Adual)
 \quad \textnormal{and} \quad
 [n]^*: \; \bbZ[x^{\pm 1}] \longrightarrow \bbZ[x^{\pm 1}]
\] 
the ring homomorphisms which are given, respectively, by the pull-back under the isogeny $[n]: \Adual \to \Adual$ and by $x^i \mapsto x^{ni}$ for $i\in \bbZ$. This is motivated by example~\ref{ex:adams}. 

\begin{thm} \label{thm:schur}
Let $P\in \Pervc(A)$ be a clean semisimple perverse sheaf, $\natural \in \{ *, \bullet\}$, and for $\natural = \bullet$ suppose the spectrum $\scrS = \scrS(P)$ is finite. Then for each partition $\alpha$ of degree $n$ we have

\[
 \rtilde_P^\natural(\alpha) \;=\; \frac{1}{n!} \sum_{\sigma \in \frakS_n} \chi_{V_\alpha}(\sigma) \cdot c_P^\natural(\sigma),
\]
with $c_P^\natural(\sigma)\in \bbZ$ given by 
\[ 
 c_P^\natural(\sigma) \;=\;
 \begin{dcases}
 \;
 \Bigl[ \, \prod_{r=1}^{\ell(\sigma)} 
 \,  
 [\sigma_r]^* (\gamma_s(P) ) \,
 \Bigr]_{s^g} 
 & \textnormal{\em for $\natural = *$}, \\ 
 \;
 \sum_{\varphi \in \scrS}
  \Bigl[ \, \iota^{-1} \prod_{r=1}^{\ell(\sigma)} 
  [\sigma_r]^* (b_x(P_\varphi)) \,
   \Bigr]_{s^g} 
 & \textnormal{\em for $\natural = \bullet$}.
 \end{dcases}
\]
In the last formula $\iota: \bbZ[s] \hookrightarrow \bbZ[x^{\pm 1}]$ denotes the embedding from section~\ref{sec:betti}.
\end{thm}

{\em Proof.} Via proposition~\ref{prop:negligible} and lemma~\ref{lem:euler}, this follows directly from lemma~\ref{lem:projector} and from the formulae for the Adams operations in example~\ref{ex:adams}. \qed

\medskip

\begin{cor} \label{cor:elliptic}
For $g=1$, any $P\in \Pervc(A)$ of Euler characteristic $\chi=\chi(P)$ and generic rank $r$ satisfies
\[
 c_P^*(\sigma) \;=\; r\, \chi^{\ell(\sigma)-1} \, \sum_{i=1}^{\ell(\sigma)} \, \sigma_i^2
 \quad \textnormal{\em and} \quad
 c_P^\bullet(\sigma) \;=\; 0
 \quad \textnormal{\em for all} \quad
 \sigma \;\in \; \frakS_n.
\]
\end{cor}

{\em Proof.} Let $p\in \rmH^2(A)$ be the class of a point. Then $\gamma(P) = \chi + r p$ by lemma~\ref{lem:euler} and the remark immediately thereafter, using our assumption that $g=1$. It follows that $c_P^*(\sigma)$ is the coefficient of $s$ in $\prod_{i=1}^{\ell(\sigma)} (\chi + \sigma_i^2 rs)$ which immediately gives the first formula. The hypercohomology of any clean perverse sheaf on an elliptic curve is concentrated in degree zero by the same argument as in the proof of lemma~\ref{lem:spectrum} below, so $b_x(P_\varphi)=\chi$ is a constant for all $\varphi \in \Pi(A)$ and hence $c_P^\bullet(\sigma)=0$. \qed

\medskip

\section{Computations for curves}

Let us illustrate the above results with the case of curves. A subvariety $Z\hookrightarrow A$ is said to {\em generate} the abelian variety $A$ if the addition map $Z^n \to A$ is surjective for all $n\gg 0$. For the perverse intersection cohomology sheaf of smooth generating curves the spectrum is trivial:

\begin{lem} \label{lem:spectrum}
If $C\hookrightarrow A$ is a smooth curve generating $A$, then $\scrS(\delta_C)=\{0\}$.  
\end{lem}

{\em Proof.} An intermediate extension argument like in the first part of the proof of proposition~\ref{prop:negligible} but with $|n|=\dim \Supp(\delta_C) = 1$ shows that any $\varphi \in \scrS(\delta_C)$ must satisfy~$\delta_{C,\varphi} \cong \delta_C$. But such an isomorphism can exist only if the character $\varphi$ is trivial, indeed the homomorphism $\pi_1(C, pt) \to \pi_1(A, 0)$ is surjective for any smooth curve $C$ that generates $A$. \qed

\medskip

Let $\theta \in \rmH^2(A)$ denote the class of an ample divisor that defines a principal polarization. We will call a smooth projective curve $C\hookrightarrow A$ a {\em Prym-Tjurin curve} of exponent $m\in \bbN$ if it generates $A$ and has fundamental class
\[ [C] \; = \; \frac{m\, \theta^{g-1}}{(g-1)!} \;\in\; \rmH^{2g-2}(A). \]
This notion comprises the case of Jacobian varieties for $m=1$ and the case of Prym varieties for $m=2$, see also~\cite[sect.~12.2]{BL}.

\begin{thm} \label{cor:prympowers}
If $C\hookrightarrow A$ is a Prym-Tjurin curve of exponent $m \in \bbN$ and $P=\delta_C$, then  
\[
 f_{P}(t) \;=\; (m^g  g!  - 1) \, t^{g-1}
 \quad \textnormal{\em and} \quad
 \ftilde_{P}(t) \;=\; m^g \, t^{g-1}  \vspace*{0.8em}
\]
and we have
\begin{align} \label{eq:schurchern}
 c^*_P(\sigma) & \;=\; 
 g!
 \, \Bigl[ \,
 \prod_{i=1}^{\ell(\sigma)} \bigl( \chi + \sigma_i^2 m s \bigr) \,
 \Bigr]_{s^g} \\ \label{eq:schurbetti}
c^\bullet_P(\sigma)  &\;=\; 
 \Bigl[ \,
 \prod_{i=1}^{\ell(\sigma)} \bigl( \chi - sq_{\sigma_i}(s) \bigr) \,
 \Bigr]_{s^g} 
\end{align}
where $\chi=\chi(P)$ and where the polynomials $q_{\sigma_i}(s)\in \bbZ[s]$ are as in definition~\ref{def:qpolynomial}.
\end{thm}

{\em Proof.} For the perverse intersection cohomology sheaf $P=\delta_C$ the spectrum is trivial by lemma~\ref{lem:spectrum}. Inserting $h_s(P) = s + \chi$ in theorems~\ref{thm:explicit} and~\ref{thm:explicit_symmetric} we therefore obtain 
\[
 f_P(t) \;=\; 
 \sum_{n=0}^g \bigl[ (\gamma_s(P)-\chi)^n  - s^n  \bigr]_{s^g} \cdot t^{g-1} (1-\chi t)^{g-n}
\]
and
\[
 \ftilde_P(t) \;=\; t^{-1} \Bigl[ \, \prod_{i=1}^g \exp \bigl( \gamma_i(P) \, tp_{2i-1}(t) \, s^i  \bigr) \Bigr]_{s^g} 
\]
where in the latter case we have used that here $\ftilde_P^\bullet(t)=0$ either by the formula in theorem~\ref{thm:explicit_symmetric} or by a direct inspection of the hypercohomology. The degree shift by one in the definition of the perverse sheaf $P=\delta_C$ gives an extra sign $-1$ in the Chern-MacPherson class $c(\delta_C) = -c(C)$. So if we denote by $p \in \rmH^{2g}(A)$ the class of a point, then
\[ 
 c(\delta_C) \; = \; \chi \cdot p \;-\; \frac{m\, \theta^{g-1}}{(g-1)!}
 \qquad \textnormal{and hence} \qquad
 \gamma_s(P) \; = \; \chi  \;+\; m\, \tdual \, s
\]
where $\tdual \in \rmH^2(\Adual)$ is the principal polarization class on $\Adual$ \cite[lemme 1]{BeauvilleFourier}. The claim about the generating polynomials is then obtained via a simple computation, and the formulae for the Schur functors follow similarly by inserting $\gamma_s(P)=\chi+m\tdual s$ and $\iota^{-1}(x^n+x^{-n})=sq_n(s)-2$ in theorem~\ref{thm:schur}.  \qed

\medskip

\bibliographystyle{amsplain}
\bibliography{Bibliography}

\providecommand{\bysame}{\leavevmode\hbox to3em{\hrulefill}\thinspace}
\providecommand{\MR}{\relax\ifhmode\unskip\space\fi MR }
\providecommand{\MRhref}[2]{%
  \href{http://www.ams.org/mathscinet-getitem?mr=#1}{#2}
}
\providecommand{\href}[2]{#2}
\begin{thebibliography}{10}

\bibitem{AtiyahPower}
{{Atiyah, M.}}, \emph{{{Power operations in $K$-theory}}}, {{Q. J. Math.}}
  \textbf{{17}} ({1966}), {165--193}.

\bibitem{BeauvillePrym}
{Beauville, A.}, \emph{{Prym varieties and the Schottky problem}}, Invent.
  Math. \textbf{41} (1977), 149--196.

\bibitem{BeauvilleFourier}
\bysame, \emph{{Quelques remarques sur la transformation de Fourier dans
  l'anneau de Chow d'une vari{\'e}t{\'e} ab{\'e}lienne}}, Lecture Notes in
  Math., vol. 1016, {Springer Verlag}, {1983}, pp.~{238--260}.

\bibitem{BL}
{Birkenhake, C. and Lange, H.}, \emph{{Complex abelian varieties (second
  ed.)}}, Grundlehren Math. Wiss., vol. 302, Springer Verlag, 2004.

\bibitem{BKdeJong}
{B\"ockle, G. and Khare, C.}, \emph{{Mod $\ell$ representations of arithmetic
  fundamental groups. II. A conjecture of A. J. de Jong.}}, {Compos. Math.}
  \textbf{{142}} ({2006}), {271--294}.

\bibitem{CappellSymmetric}
{{Cappell, S. E., Maxim, L., Sch{\"u}rmann, J., Shaneson, J. L. and Yokura,
  S.}}, \emph{{Characteristic classes of symmetric products of complex
  quasi-projective varieties}}, {J. Reine Angew. Mathematik} ({2015}),
  {\url{doi:10.1515/crelle-2014-0114}}.

\bibitem{DebTorelli}
{Debarre, O.}, \emph{{Sur la d{\'e}monstration de A.~Weil du th{\'e}or{\`e}me
  de Torelli pour les courbes}}, Compos. Math. \textbf{58} (1986), 3--11.

\bibitem{DebarreSchottkyUpdate}
\bysame, \emph{{The Schottky problem --- an update}}, Complex algebraic
  geometry, Math. Sci. Res. Inst. Publ., vol.~28, 1995, pp.~57--64.

\bibitem{DelCT}
{Deligne, P.}, \emph{{Cat{\'e}gories tensorielles}}, Moscow Math. J. \textbf{2}
  (2002), 227--248.

\bibitem{DonagiPrym}
{Donagi, R.}, \emph{{The fibers of the Prym map}}, Contemporary Mathematics
  \textbf{136} (1992), 55--125.

\bibitem{DrinfeldKashiwara}
{Drinfeld, V.}, \emph{{On a conjecture of Kashiwara}}, Math. Res. Lett.
  \textbf{8} (2001), 713--728.

\bibitem{GaL}
{Gabber, O. and Loeser, F.}, \emph{{Faisceaux pervers $\ell$-adiques sur un
  tore}}, Duke Math. J. \textbf{83} (1996), 501--606.

\bibitem{GaitsgoryDeJong}
{Gaitsgory, D.}, \emph{{On de Jong's conjecture}}, Israel J. Math. \textbf{157}
  (2007), 155--191.

\bibitem{Heinloth}
{Heinloth, F.}, \emph{{A note on functional equations for zeta functions with
  values in Chow motives}}, {Ann. Inst. Fourier (Grenoble)} \textbf{{57}}
  ({2007}), {1927--1945}.

\bibitem{HirzebruchEulerian}
{Hirzebruch, F.}, \emph{{Eulerian polynomials}}, {M{\"u}nster J. of Math.~1}
  ({2008}), 9--14.

\bibitem{IzadiGeometricStructure}
{Izadi, E.}, \emph{{The geometric structure of $\mathcal{A}_4$, the structure
  of the Prym map, double solids and $\Gamma_{00}$-divisors}}, J.~Reine
  Angew.~Math. \textbf{462} (1995), 93--158.

\bibitem{KatzSatoTate}
{Katz, N. M.}, \emph{{Convolution and equidistribution: Sato-Tate theorems for
  finite-field Mellin transforms}}, {Ann. of Math. Stud.}, vol. {180},
  {Princeton Univ. Press}, {2012}.

\bibitem{KnutsonLambda}
{{Knutson, D.}}, \emph{{{$\lambda$-rings and the representation theory of the
  symmetric group}}}, {{Lecture Notes in Math.}}, vol. {308}, {Springer
  Verlag}, {1973}.

\bibitem{KraemerSemiabelian}
{Kr{\"a}mer, T.}, \emph{{Perverse sheaves on semiabelian varieties}}, {Rend.
  Semin. Mat. Univ. Padova} \textbf{{132}} ({2014}), {83--102}.

\bibitem{KrE6}
\bysame, \emph{{On a family of surfaces of general type attached to abelian
  fourfolds and the Weyl group $W(E_6)$}}, {Int. Math. Res. Notices IMRN}
  \textbf{{2015, 24}} ({2015}), {13062--13105}.

\bibitem{KrCubic}
{{Kr{\"a}mer, T.}}, \emph{{Cubic threefolds, Fano surfaces and the monodromy of
  the Gauss map}}, {Manuscripta Math.} \textbf{{149}} ({2016}), {303--314}.

\bibitem{KrWSchottky}
{Kr{\"a}mer, T. and Weissauer, R.}, \emph{{The Tannaka group of the theta
  divisor on a generic principally polarized abelian variety}}, {Math. Z.}
  \textbf{{281}} ({2015}), {723--745}.

\bibitem{KrWVanishing}
\bysame, \emph{{Vanishing theorems for constructible sheaves on abelian
  varieties}}, {J. Alg. Geom.} \textbf{{24}} ({2015}), {531--568}.

\bibitem{Kwiecinski_Formule}
{Kwieci\'nski, M.}, \emph{{Formule du produit pour les classes
  caract\'eristiques de Chern-Schwartz-MacPherson et homologie
  d'intersection}}, {C. R. Acad. Sci. Paris} \textbf{{314}} ({1992}),
  {625--628}.

\bibitem{KY_Product}
{Kwieci\'nski, M. and Yokura, S.}, \emph{{Product formula for twisted
  MacPherson classes}}, {Proc. Japan Acad. Ser. A Math. Sci.} \textbf{{68}}
  ({1992}), {167--171}.

\bibitem{MacPhersonChern}
{MacPherson, R. D.}, \emph{{Chern classes for singular algebraic varieties}},
  {Ann. of Math.} \textbf{{100}} ({1974}), {423--432}.

\bibitem{MaximSchuermannTwisted}
{Maxim, L. and Sch\"urmann, J.}, \emph{{Twisted genera of symmetric products}},
  {Selecta Math. (N.S.)} \textbf{{18}} ({2012}), {283--317}.

\bibitem{MochizukiAsymptotic}
{Mochizuki, T.}, \emph{{Asymptotic behaviour of tame harmonic bundles and an
  application to pure twistor $D$-modules}}, Mem. Amer. Math. Soc. \textbf{185}
  (2007).

\bibitem{SabbahPolarizable}
{Sabbah, C.}, \emph{{Polarizable Twistor $\mathscr{D}$-modules}},
  {Ast{\'e}risque} \textbf{300} (2005).

\bibitem{SchnellHolonomic}
{Schnell, C.}, \emph{{Holonomic $\mathscr{D}$-modules on abelian varieties}},
  {Publ. Math. Inst. Hautes {\'E}tudes Sci.} \textbf{{121}} ({2015}), {1--55}.

\bibitem{SchuermannTopology}
{Sch\"urmann, J.}, \emph{{Topology of singular spaces and constructible
  sheaves}}, {Monografie Matematyczne}, vol.~{63}, {Birkh\"auser}, {2003}.

\bibitem{WeBN}
{Weissauer, R.}, \emph{{Brill-Noether sheaves}}, {\url{arXiv:math/0610923}}.

\bibitem{WeTorelli}
\bysame, \emph{{Torelli's theorem from the topological point of view}}, Modular
  forms on Schiermonnikoog ({Edixhoven, B., van der Geer, G. and Moonen, B.},
  ed.), Cambridge University Press, 2008, pp.~275--284.

\end{thebibliography}

\end{document}